\documentclass[12pt,a4paper,reqno]{amsart}

\usepackage{amssymb}
\usepackage{latexsym}
\usepackage{exscale}

\numberwithin{equation}{section}
\newtheorem{theor}{Theorem}[section]
\newtheorem{lemma}[theor]{Lemma}

\newtheorem{remark}[theor]{Remark}
\newtheorem{prop}[theor]{Proposition}

\newcommand{\ZZ}{\mathbb{Z}}
\newcommand{\RR}{\mathbb{R}}
\newcommand{\NN}{\mathbb{N}}

\newcommand{\W}{\mathcal{W}}

\newcommand{\T}{\mathcal{T}}
\newcommand{\U}{\mathcal{U}}
\newcommand{\V}{\mathcal{V}}
\newcommand{\WW}{\mathbb{W}}
\newcommand{\UU}{\mathbb{U}}
\newcommand{\VV}{\mathbb{V}}
\newcommand{\TT}{\mathbb{T}}
\newcommand{\KK}{\mathbb{K}}

\newcommand{\D}{\mathcal{D}}

\newcommand{\ep}{\varepsilon}

\renewcommand{\emptyset}{\mbox{\rm \O}}

\begin{document}
\allowdisplaybreaks

\title{BCR algorithm and the $T(b)$ theorem}
      \thanks{Research supported by
NNSF of China No.10001027, the innovation funds of Wuhan
University and the subject construction funds of Mathematic and
Statistic School, Wuhan University.}
\author {P. Auscher}

\address{Pascal Auscher
\\
Universit\'e de Paris-Sud, UMR du CNRS 8628,
\\
91405 Orsay Cedex, France} \email{pascal.auscher@math.u-psud.fr}
 
\author {Q.X.Yang}
\address{Qi Xiang Yang 
\\ Mathematic and statistic school, Wuhan University, 430072 Hubei, China}
\email{qxyang@whu.edu.cn}

\date{\today}


\begin{abstract}  We show using  the Beylkin-Coifman-Rokhlin algorithm in the Haar basis that  any singular integral operator can be written as the sum of a bounded operator on $L^p$, $1<p<\infty$. and of a perfect dyadic singular integral operator.  This allows to deduce a  local $T(b)$ theorem for singular integral operators from the one  for perfect dyadic singular integral operators obtained by Hofmann, Muscalu,  Thiele, Tao and the first author. \end{abstract}

\maketitle

\section{Introduction}

The purpose of this note is to  fill in a gap of \cite{AHMTT} concerning  a local $T(b)$ theorem for singular integrals with a method that could be of interest elsewhere.

In \cite{C}, M. Christ  proves a local $T(b)$ theorem for singular integral operators on a space of homogeneous type, the motivation being the potential application  to several questions related to analytic capacity. It lead to the solution of  the Vitushkin's conjecture by G. David \cite{D} or  to a proof of the semiadditivity
of analytic capacity (Painlev\'e problem)   by X. Tolsa \cite{T}.   Those solutions required similar $T(b)$ theorems but  in non-homogeneous spaces  as developed by G. David \cite{D}, and S. Nazarov, S. Treil and A. Volberg \cite{NTV1, NTV2, V}.

Let us explain Christ's theorem and the word ``local''. He introduces the notion of accretive systems $(b_Q)$ consisting   of  functions supported on the corresponding cube $Q$, bounded,  non-degenerate  (i.e. of mean 1). He requires that each $b_Q$ is mapped through the operator  to a bounded function  on $Q$ (and a similar hypothesis for the adjoint with a different accretive system if need be).     He designs globally defined  para-accretive functions $b$ and $b^*$ adapted to the operator  and its adjoint,  and  applies the 
David-Journ\'e-Semmes' $T(b)$ theorem \cite{DJS} to obtain the $L^2$ boundedness of the operator. 

A generalization of Christ's result is proposed in \cite{AHMTT} in Euclidean space for a model situation.

\begin{theor}[\cite{AHMTT}]\label{theor:ahmttintro} Assume that $T$ is a  perfect dyadic singular integral operator. Assume that there exists a constant $C\ge 0$ such that for 
each dyadic cube $Q$, one can find functions $b^{1}_{Q}, b^{2}_{Q}$ supported in $Q$ with 
\begin{equation}
 \int_Q  b^{1}_{Q} = |Q|= \int_Q  b^{2}_{Q}, 
\end{equation}
\begin{equation}\label{eq:bQintro}
\int_{Q} |b^{1}_{Q}|^{2}
+|b^{2}_{Q}|^{2}  \leq C|Q|
\end{equation}
\begin{equation}\label{eq:TbQintro}
\int_{Q}  |Tb^{1}_{Q}|^{2} +|T^{*}b^{2}_{Q}|^{2}
\leq C|Q| .
\end{equation}
Then $T$ is bounded on $L^2(\RR^n)$. 
\end{theor}
The idea of proof  is different from Christ's argument (in fact, it is not clear how to adapt it): it amounts to verify the hypotheses of  a variant of the $T(1)$ theorem of David-Journ\'e \cite{DJ}, namely a local $T(1)$ theorem.  Perfect dyadic means essentially that the regularity is adapted to the dyadic grid: any function supported in a dyadic cube with mean 0 is mapped to a function supported in the same cube. This property kills most tail terms that  would appear with standard singular integrals. 

The following natural extension  is announced in \cite{AHMTT}.

\begin{theor}\label{theor:tblocdyadicintro} Assume that $T$ is singular integral operator with locally bounded  kernel on $\RR^n\times \RR^n$. Assume that there exists a constant $C\ge 0$ such that for 
each dyadic cube $Q$, one can find functions $b^{1}_{Q}, b^{2}_{Q}$ supported on $Q$ with
\begin{equation}
 \int_Q  b^{1}_{Q} = |Q|= \int_Q  b^{2}_{Q}, 
\end{equation}
\begin{equation}\label{eq:bQbisintro}
\int_{Q} |b^{1}_{Q}|^{2}
+|b^{2}_{Q}|^{2}  \leq C|Q|
\end{equation}
\begin{equation}\label{eq:TbQbisintro}
\int_{Q}  |Tb^{1}_{Q}|^{2} +|T^{*}b^{2}_{Q}|^{2}
\leq C|Q| .
\end{equation}
Then $T$ is bounded on $L^2(\RR^n)$. 
\end{theor}

It looks like a straightforward  exercise to adapt the proof in the model case by handling the tails as error terms. This is actually said   in \cite{AHMTT} but, on second thoughts,  it may have been too optimistic\footnote{The first author feels responsible for that.}. The far away tails are indeed easy to handle,  that is integrals $\int_{Q\times R} g(x)K(x,y) f(y) \, dxdy$ when $R \cap 3Q\ne \emptyset$ or $Q \cap 3R\ne \emptyset$ with $f$ or $g$ having mean value 0. But the same integrals on adjacents cubes of different sizes  seem a problem. The reader can be convinced by reading the proof of   Theorem 8.6  of \cite{AAAHK} in \cite{H} where  the hypothesis \eqref{eq:bQbisintro} has been  strengthened  to work out the transposition of the \cite{AHMTT} argument. 

It would be interesting to have a direct proof of this result but we have not succeeded.  
Our idea is to  reduce to the model case via the following result, interesting on its own. 

\begin{theor}\label{theor:siotopdsio:intro} Assume that $T$ is a singular integral operator with  locally bounded kernel on $\RR^n\times \RR^n$. Then there exists a  perfect dyadic singular integral operator
 $\TT$    such that $T-\TT$ is  bounded on $L^p(\RR^n)$ for all  $1<p<\infty$.
\end{theor}

 This is done using the Beylkin-Coifman-Rokhlin algorithm in the Haar basis and ideas from the PhD thesis of one of us \cite{Y1}.

Let us say that the extension of Christ's result for singular integrals is not just an academic exercice.  Such a generalization found recently an application in \cite{AAAHK}  towards the $L^2$ boundedness of boundary layer potentials for some PDE's. Other potential applications require a similar theorem where $L^2$ conditions on the accretive system  replaced by $L^p$ conditions for $p>1$\footnote{Personal communication of Steve Hofmann.}. For perfect dyadic models, it is  remarked in \cite{AHMTT} that the $L^2$ conditions  can be replaced  by $L^p$ conditions  for any  $1<p\le \infty$\footnote{The inequality $1\le p$ is written but this is obviously a typo as the whole argument depends on the stopping time argument in Lemma 6.5, which does not give anything for $p=1$.}.  See the extension of Theorem \ref{theor:tblocdyadicintro} in the text. At present, none of the arguments for standard singular integrals  in   \cite{H} or here work  with $L^p$ conditions for $p<2$. We leave this question  open.

\section{From an operator to a perfect dyadic operator}

Here is a formal approach. We begin with the BCR algorithm in the Haar basis. Consider the Haar wavelets in $\RR^n$ given by
\begin{equation}
\psi_{j,k}(x)= 2^{nj/2}\psi(2^jx-k), \ j\in \ZZ, k\in \ZZ^n, \psi \in E
\end{equation}
where $E$ is a set of cardinal $2^n-1$. Recall that $\psi_{j,k}$ has support in the dyadic cube $Q=Q_{j,k}= 2^{-j}k +2^{-j}[0,1)^{n}$, that $\int \psi_{j,k}=0$ and that $\{\psi_{j,k}\}$ is an orthonormal basis of $L^2(\RR^n)$.  Define also 
$$\phi_{j,k}(x)= 2^{nj/2}\phi(2^jx-k), \ j\in \ZZ, k\in \ZZ^n, \phi=1_{[0,1)^n}.
$$
We also use the notation $\psi_Q$ and $\phi_Q$ when more convenient.  It is understood that  that the Haar functions $\psi$ describe the set $E$ and we  forget from now on to mention this as it plays no role.  

For $j\in \ZZ$, we let $V_j$ be the closed subspace of $L^2$ generated by the orthonormal system $\phi_{j,k}$, $k\in \ZZ^n$ and $W_j$ the closed subspace of $L^2$ generated by the orthonormal system $\psi_{j,k}$, $k\in \ZZ^n$.   It is well-known that $V_j$ and $W_j $ are orthogonal spaces and $L^2(\RR^n)= \oplus W_j$. Furthermore, one has  $P_{j+1}=P_j + Q_j$   where  $P_j$ and $Q_j$ are the orthogonal projections onto $V_j$ and $W_j$. In what follows, $\langle \ , \  \rangle$ denotes the  bilinear duality bracket and the adjoint of an operator $T$ for this duality is denoted by $T^*$.

Consider  an operator   $T$ for which one can define the coefficients for all $j\in \ZZ$, 
\begin{equation}\label{eq:d}
\langle \phi_Q, T\phi_R\rangle, \quad Q=Q_{j,k}, R= Q_{j,\ell}, 
\end{equation}
\begin{equation}\label{eq:a}
\langle \psi_Q, T\psi_R\rangle = a_{Q,R}= a^j_{k,\ell}, \quad Q=Q_{j,k}, R= Q_{j,\ell}, 
\end{equation}
\begin{equation}\label{eq:b}
\langle \psi_Q, T\phi_R \rangle = b_{Q,R}= b^j_{k,\ell}, \quad Q=Q_{j,k}, R= Q_{j,\ell}, 
\end{equation}
\begin{equation}\label{eq:c}
\langle \phi_Q, T\psi_R \rangle = c_{Q,R}= c^j_{k,\ell}, \quad Q=Q_{j,k}, R= Q_{j,\ell}, 
\end{equation}
and such that  for $f,g$ in some appropriate  vector space(s) of measurable functions, 
\begin{equation}\label{eq:lim+}
\lim_{j\to +\infty}\langle P_j g, T P_j f \rangle = \langle g, T f \rangle \end{equation}
and 
\begin{equation}\label{eq:lim-}
\lim_{j\to -\infty}\langle P_j g, T P_j f \rangle =0. \end{equation}
Note that $\langle P_j g, T P_j f \rangle$ is defined using the first set of coefficients in \eqref{eq:d}.  
Then, one can expand formally
\begin{equation}
 \langle g, T f \rangle =   \langle g, U f \rangle + \langle g, V f \rangle + \langle g, W f \rangle
  \end{equation}
where
\begin{equation}
 \langle g, U f \rangle =   \sum_{j=-\infty}^{+\infty}  \langle Q_j g, T Q_j f \rangle
 \end{equation}
\begin{equation}
 \langle g, V f \rangle =   \sum_{j=-\infty}^{+\infty}  \langle Q_j g, T P_j f \rangle
 \end{equation}
\begin{equation}
 \langle g, W f \rangle =   \sum_{j=-\infty}^{+\infty}  \langle P_j g, T Q_j f \rangle.
 \end{equation}
 Expanding on the bases of $V_j$ and $W_j$, one finds
 \begin{equation}
 \langle g, U f \rangle =   \sum_{j=-\infty}^{+\infty}\sum_{k,\ell\in \ZZ^n}  \langle  g,  \psi_{j,k} \rangle  a^j_{k,\ell}   \langle   \psi_{j,\ell}, f \rangle
 \end{equation}
\begin{equation}
 \langle g, V f \rangle =   \sum_{j=-\infty}^{+\infty}\sum_{k,\ell\in \ZZ^n}  \langle  g,  \psi_{j,k} \rangle  b^j_{k,\ell}   \langle   \phi_{j,\ell},  f \rangle
 \end{equation}
 \begin{equation}
 \langle g, W f \rangle =   \sum_{j=-\infty}^{+\infty}\sum_{k,\ell\in \ZZ^n}  \langle  g,  \phi_{j,k} \rangle  c^j_{k,\ell}   \langle   \psi_{j,\ell} , f \rangle.
 \end{equation}
 This is the so-called BCR algorithm in the Haar basis. The operator $U$ is diagonal in the decomposition of $L^2$ given by the $W_j$. The operator $V$ is some sort of paraproduct and $W$ is like $V^*$.  This decomposition can be used to prove the T(1) theorem. 
 
 Let us go further and modify formally  $U$, $V$, $W$. 
 Set
\begin{equation}
\alpha^j_{k,\ell} =  \begin{cases}  a^j_{k,\ell}\ , & \mathrm{if}\   k\ne \ell,
\\ 0\ ,  & \mathrm{if}\  k=\ell.
\end{cases}
\end{equation}
 \begin{equation}
\beta^j_{k,\ell} =  \begin{cases}  b^j_{k,\ell}\ , & \mathrm{if}\   k\ne \ell,
\\ - \sum_{m\ne 0} b^j_{k,k+m}\ ,  & \mathrm{if}\  k=\ell.
\end{cases}
\end{equation}
 \begin{equation}
\gamma^j_{k,\ell} =  \begin{cases}  c^j_{k,\ell}\ , & \mathrm{if}\   k\ne \ell,
\\ - \sum_{m\ne 0} c^j_{k+m,k}\ ,  & \mathrm{if}\  k=\ell.
\end{cases}
\end{equation}
and $\U, \V, \W$ the operators associated with the family of coefficients $\alpha, \beta, \gamma$ as $U,V,W$ with the family of  coefficients $a,b,c$. The $\alpha, \beta, \gamma$ are designed so that $\V(1)=\V^*(1)=0$ and $\W(1)=\W^*(1)=0$ (of course, this has only a formal meaning) so that  the main result  is the following.

\begin{theor}\label{theor:bcmr}
Assume that for some $s>0$ and $C>0$ one has for all $j$, $k,\ell$ with $k\ne \ell$, 
\begin{equation}\label{eq:bcmr}
|a^j_{k,\ell}| + |b^j_{k,\ell}| + |c^j_{k,\ell}| \le C (1+|k-\ell |)^{-n-s}.
\end{equation}
Then $\U, \V,\W$ are bounded operators on $L^p(\RR^n)$, $1<p<\infty$ and also from $H^1_d
(\RR^n)$ into $L^1(\RR^n)$ and from $L^\infty(\RR^n) $ into BMO$_d(\RR^n)$. \end{theor}

We set $\T=\U+\V+\W$. Here, $H^1_d(\RR^n)$ and BMO$_d(\RR^n)$ are the dyadic Hardy and BMO space.  The proof is in section \ref{sec:proof}. In fact, a decay in $|k-\ell|^{-n} \ln^{-2-\ep} (1+ |k-\ell|)$  with $\ep>0$ suffices.
\

We remark that the point of this statement is to avoid use of the ``diagonal coefficients'' in the families $a,b,c$ as this would require some sort of weak boundedness property which we do not want to assume. 

\

This theorem has its origin in \cite{Y1}  where the Haar functions are replaced by smooth compactly supported wavelets.   But 
the point was different.  The operator $T$ was supposed bounded on $L^2$ and the objective was  to obtain the rate of approximation of $T$ by some truncated $T_m$ in  the non-standard representation defined by the BCR algorithm. Here, we do not assume that the original $T$ is bounded. 
See also \cite{Y2,Y3, DYY} for related ideas. 

\

Let $\UU, \VV, \WW$ be the differences $U-\U, V-\V, W-\W$ and  $\TT=\UU+\VV+\WW$.  
Thus the boundedness of $T$ on $L^2$ is equivalent to that of $\TT$. Note that \begin{equation}
 \langle g, \UU f \rangle =   \sum_{j=-\infty}^{+\infty}\sum_{k\in \ZZ^n}  \langle  g,  \psi_{j,k} \rangle  \mathbf{a} ^j_{k}   \langle   \psi_{j,k}, f \rangle
 \end{equation}
\begin{equation}
 \langle g, \VV f \rangle =   \sum_{j=-\infty}^{+\infty}\sum_{k\in \ZZ^n}  \langle  g,  \psi_{j,k} \rangle  \mathbf{b} ^j_{k}   \langle   \phi_{j,k},  f \rangle
 \end{equation}
 \begin{equation}
 \langle g, \WW f \rangle =   \sum_{j=-\infty}^{+\infty}\sum_{k\in \ZZ^n}  \langle  g,  \phi_{j,k} \rangle \mathbf{c} ^j_{k}   \langle   \psi_{j,k} , f \rangle
 \end{equation}
 for some family of complex coefficients $\mathbf{a}, \mathbf{b}, \mathbf{c}$. The only use of these formulae is in the following (formal) observation.
 
 \begin{lemma}\label{lemma:perfect}
If $f$ is supported in a dyadic cube and has mean 0, then $\TT f$ is supported in the same cube in the sense that $\langle g, \TT f\rangle=0$ if $g$ is supported away from $Q$. 
\end{lemma}

Let $Q$ be the dyadic cube supporting $f$. The coefficients
$ \langle   \psi_{j,k} , f \rangle$ and $ \langle   \phi_{j,k} , f \rangle$ are 0 if $Q_{j,k} \cap Q=\emptyset$  and also if $Q \subsetneq Q_{j,k}$ since $f$ has mean $0$. Hence the sums reduce to couples $(j,k)$ such that $Q_{j,k} \subset Q$. Thus, if $g$ is supported away from $Q$, we have $ \langle g, \TT f \rangle =0$. 

\

We are now ready to apply all this to singular integral operators.

\section{Application to singular integral operators}

Assume that $T$ is a singular integral operator, that is a linear continuous operator  from $\D(\RR^n)$ to $\D'(\RR^n)$ whose distributional kernel $K(x,y)$ satisfies the  Calder\'on-Zygmund estimates,  that  is 
the size condition
\begin{equation}\label{eq:size}
|K(x,y)| \leq C|x-y|^{-n}, 
\end{equation}
for all $x,y$ with $x\ne y$ and the   regularity condition for some $0<s<1$
\begin{equation}\label{eq:regularity}
|K(x,y)-K(x',y)| +|K(y,x)-K(y,x')| \leq
C\frac{|x-x'|^{s}}{ |x-y|^{n+s}},\end{equation}
for all $x,x',y $ {with} $|x-x'|\leq
\frac{1}{2}|x-y|.$

Assume also that $K$ is  locally bounded on $\RR^n\times \RR^n$. The local boundedness of $K$ guarantees that one can start the BCR algorithm with $T$ and obtain  operators $\T$ and $\TT$.  More precisely, we  first 
extend $\langle g, Tf \rangle$  a priori defined for $f,g \in \D(\RR^n)$  to $f,g \in L^1_c(\RR^n)$, the space of compactly supported integrable functions, by 
$$
\langle g, Tf \rangle = \iint_{\RR^n\times \RR^n} g(x) {}{K(x,y) f(y)} \, dxdy.
$$
Hence all the coefficients $a$, $b$, $c$   in \eqref{eq:a}, \eqref{eq:b}, \eqref{eq:c} can be computed and the limits in \eqref{eq:lim+} and \eqref{eq:lim-} hold for $f,g \in L^1_c(\RR^n)$.
Moreover, the Calder\'on-Zygmund conditions on the kernel and standard estimates insure that \eqref{eq:bcmr} hold so that Theorem \ref{theor:bcmr} applies. Thus, the operator $\T$ is bounded on $L^p(\RR^n)$, $1<p<\infty$. Furthermore, one has 

\begin{prop}
The distribution kernel  of $\T$ satisfies the size condition
 \eqref{eq:size}.
\end{prop}

This is also a standard computation from \eqref{eq:bcmr}. Hence by difference and incorporating lemma \ref{lemma:perfect},  $\TT$ has the following properties:
\begin{enumerate}
\item $\TT$ is a linear continuous operator  from $\D(\RR^n)$ to $\D'(\RR^n)$. 
\item $\TT$ has a kernel satisfying the size condition \eqref{eq:size}.
\item $\langle g, \TT f\rangle$ is well-defined for pairs of functions $(f,g) \in  L^p_c(\RR^n) \times L^{p'}_c(\RR^n)$ for $1<p<\infty$  and  if, furthermore, they are integrable with  support on disjoint cubes (up to a set of measure 0)
$$\langle g, \TT f\rangle= \iint_{\RR^n\times \RR^n} g(x) {}{\KK(x,y) f(y)} \, dxdy.
$$
\item for all $(f,g)$ as above, if $f$ has support in a dyadic cube and mean 0, then 
$\langle g, \TT f\rangle=0$ when the support of $g$ does not meet $Q$ (up to a set of measure 0).
\end{enumerate}

We say that an operator satisfying the above four properties is a  perfect dyadic singular integral operator. We note that this is not exactly the definition in \cite{AHMTT}, which is concerned with a dyadic and finite model, where the operator is defined on a finite dimensional subspace of the one generated by the $\psi_Q$ and the $\phi_Q$. But this is a superficial difference.
Let us summarize the main result. 

\begin{theor}\label{theor:siotopdsio} Assume that $T$ is a singular integral operator with  locally bounded kernel on $\RR^n\times \RR^n$. Then there exists a  perfect dyadic singular integral operator
 $\TT$    such that $T-\TT$ is  bounded on $L^p(\RR^n)$ for all  $1<p<\infty$.
\end{theor}

The  criterion for $L^2$ boundedness of perfect dyadic singular integral operators in \cite{AHMTT} is (See  Theorem 6.8 there  when $p=q=2$  and  a remark after the proof  for general $p,q$).

\begin{theor}\label{theor:ahmtt} Assume that $T$ is a  perfect dyadic singular integral operator. Let 
$1<p,q\leq\infty$ with dual exponents $p',q'$. Assume that there exists a constant $C\ge 0$ such that for 
each dyadic cube $Q$, one can find functions $b^{1}_{Q}, b^{2}_{Q}$ supported in $Q$ with 
\begin{equation}
 \int_Q  b^{1}_{Q} = |Q|= \int_Q  b^{2}_{Q}, 
\end{equation}
\begin{equation}\label{eq:bQ}
\int_{Q} |b^{1}_{Q}|^{p}
+|b^{2}_{Q}|^{q}  \leq C|Q|
\end{equation}
\begin{equation}\label{eq:TbQ}
\int_{Q}  |Tb^{1}_{Q}|^{q'} +|T^{*}b^{2}_{Q}|^{p'}
\leq C|Q| .
\end{equation}
Then $T$ is bounded on $L^2(\RR^n)$. 
\end{theor}

Although we have a different definition of  perfect dyadic operators, the  proof there can be copied \textit{in extenso} in our case. 
The non trivial part is to prove first
\begin{equation}\label{eq:T1}
\int_{Q}  |T\mathbf{1}_{Q}| +|T^{*}\mathbf{1}_{Q}|
\leq C'|Q| .
\end{equation}
Then, one deduces $L^2$ boundedness by a version of the $T(1)$ theorem for dyadic perfect operators.

The conclusion of this discussion is the following local $T(b)$ theorem for singular integral operators. 

\begin{theor}\label{theor:tblocdyadic} Assume that $T$ is a  perfect dyadic singular integral operator. Let 
$1<p,q\leq\infty$ with dual exponents $p',q'$ be such that $1/p + 1/q\leq 1$. Assume that there exists a constant $C\ge 0$ such that for 
each dyadic cube $Q$, one can find functions $b^{1}_{Q}, b^{2}_{Q}$ supported in $Q$ with 
\begin{equation}
 \int_Q  b^{1}_{Q} = |Q|= \int_Q  b^{2}_{Q}, 
\end{equation}
\begin{equation}\label{eq:bQ}
\int_{Q} |b^{1}_{Q}|^{p}
+|b^{2}_{Q}|^{q}  \leq C|Q|
\end{equation}
\begin{equation}\label{eq:TbQ}
\int_{Q}  |Tb^{1}_{Q}|^{q'} +|T^{*}b^{2}_{Q}|^{p'}
\leq C|Q| .
\end{equation}
Then $T$ is bounded on $L^2(\RR^n)$. 
\end{theor}

Here is the proof.  Write $T=\T+\TT$. Since $\T$ is bounded on $L^p$ and $q'\le p$, by \eqref{eq:bQ} we have
$$
\bigg(\frac{1}{|Q|} \int_{Q}  |\T b^{1}_{Q}|^{q'} \bigg)^{1/q'} \le   \bigg(\frac{1}{|Q|}\int_{Q}  |\T b^{1}_{Q}|^{p} \bigg)^{1/p} \le 
\|\T\|_{p,p}\bigg(\frac{1}{|Q|} \int_{Q}  | b^{1}_{Q}|^{p} \bigg)^{1/p} \leq C \|\T\|_{p,p}.
$$
Thus the same conclusion holds for $\TT b^1_Q$ by \eqref{eq:TbQ} with constant
$C \|\T\|_{p,p} +C$. Similarly
$$
\bigg(\frac{1}{|Q|} \int_{Q}  |\TT^* b^{2}_{Q}|^{p'} \bigg)^{1/p'}  \leq C\|\T^*\|_{q,q} +C  .
$$
Hence we can apply  Theorem \ref{theor:tblocdyadic} to $\TT$ and conclude that 
$\TT$, hence $T$,  is bounded on $L^2(\RR^n)$.

\begin{remark}
We do not know how to  drop  the constraint $1/p+1/q\leq 1$. It is satisfied if $p=q=2$, which proves  Theorem \ref{theor:tblocdyadicintro}. \end{remark} 

\begin{remark}
If one does not want to develop the $T(1)$ theory for perfect dyadic operators, here is a direct way: first, prove \eqref{eq:T1} for $\TT$ following \cite{AHMTT}, then observe that this yields back the same conclusion for $T$. This classically implies the $L^2$ boundedness of $T$ by the $T(1)$ theorem for singular integral operators.   \end{remark} 

\begin{remark}
Actually the Calder\'on-Zygmund conditions on the kernel of $T$ can be weakened. It suffices that  for $Q,R$ distinct dyadic cubes with same sizes
\begin{equation}\label{eq:weaksize}
\int_Q \int_R |K(x,y)| \, dx dy \le C |Q|
\end{equation}
whenever   $Q$ and $R$ are adjacent $($ie $d(Q,R)=0)$ and
$$
\int_Q\int_R |K(x,y)-K(x,y_R)| \, dxdy \le \frac C   {{d(Q,R)} ^{n}} \ln^{-2-\ep} \bigg(2+ \frac {d(Q,R)} {|Q|^{1/n}}\bigg)
$$
where $y_R$ is the center of $R$, otherwise $($i.e., $d(Q,R)>0)$, and similarly for $K(y,x)$. It is easy to adapt Theorem 
\ref{theor:siotopdsio} with such hypotheses.  In such a case the kernel of $\T$ satisfies \eqref{eq:size} and the kernel of $\TT$ \eqref{eq:weaksize}. Next, 
the proof of  Theorem \ref{theor:ahmtt} easily  adapts under  \eqref{eq:weaksize} by changing the conclusion of Corollary 6.10 in \cite{AHMTT} to, with the notation there, $|\langle T(b^1_P \chi_{I_Q}), \chi_{2I_Q}\rangle |\lesssim K |I_Q|$, as this suffices to run the argument.  \end{remark}

\section{proof of Theorem \ref{theor:bcmr}}\label{sec:proof}

The case of $\U$ is the easiest one. In fact, it is bounded on all $L^p$, $1<p<\infty$,  on
$H^1_d$  and on BMO$_d$.  This is classical but we include a proof for convenience. Let us see the $L^2$ boundedness first. Set
$$
A= \sup_{j, k} \bigg\{\sum_{\ell } |\alpha^{j}_{k,\ell}| + |\alpha^{j}_{\ell,k}|\bigg\}.
$$
Recall that $\alpha^j_{k,k}=0$ so that by \eqref{eq:bcmr}, $A<\infty$.
 Write $f =  \sum_{j}f_{j}$ with $f_j=Q_j f$. Then, by Schur's lemma and 
 using the orthonormal basis property of the Haar functions,
  $$
\bigg| \sum_{k,\ell\in \ZZ^n}  \langle  g,  \psi_{j,k} \rangle  \alpha^{j}_{k,\ell}   \langle   \psi_{j,\ell},  f \rangle\bigg| \le A \|g_j\|_2\|f_j\|_2.
$$
Hence
$$
| \langle g, \U f \rangle|  \le A \sum_{j=-\infty}^{\infty}  \|g_j\|_2\|f_j\|_2 \le A \|g\|_2 \|f\|_2.
$$

\

It remains to prove the $H^1_d$ boundedness of $\U$ as the boundedness on BMO$_d$ is obtained by duality and the $L^p$ boundedness by interpolation. To do that, we pick an $L^2$ dyadic atom $a$: it is supported in a dyadic cube $Q$, its $L^2$ norm is bounded  by $1/|Q|^{1/2}$ and it is of mean 0. By scale and translation invariance, it suffices to assume that $Q=Q_{0,0}$.  Write $a=\sum_{Q_{j,\ell} \subset Q_{0,0}} \langle a, \psi_{j,\ell}\rangle  \psi_{j,\ell}$ so that  $\|a\|_2^2 = \sum_{Q_{j,\ell} \subset Q_{0,0}} |\langle a, \psi_{j,\ell}\rangle|^2 $. We have
 \begin{align*}
 \U a&= \sum_{j=-\infty}^\infty \sum_{k,\ell}   \langle a, \psi_{j,\ell} \rangle {}{\alpha^{j}_{k,\ell} } \psi_{j, k} 
 \\
 &=  \sum_{j=0}^\infty \sum_{k}\sum_{\ell; Q_{j,\ell} \subset Q_{0,0}} \langle a, \psi_{j,\ell} \rangle {}{\alpha^{j}_{k,\ell}} \psi_{j, k}
 \\
 &
 = \sum_{m\in \ZZ^n}a_m
 \end{align*}
 with 
 $$
 a_m =   \sum_{j=0}^\infty  \sum_{k \, ; \, Q_{j,k} \subset Q_{0,m}}\bigg\{ \sum_{\ell; Q_{j,\ell} \subset Q_{0,0}} \langle a, \psi_{j,\ell} \rangle {}{\alpha^{j}_{k,\ell}}\bigg\} \psi_{j, k}.
 $$
 We have that $a_m$ is supported in $Q_{0,m}$ and has mean 0. Thus  $\|a_m\|^{-1}_2 a_m$ is an $L^2$ dyadic atom. It suffices to show that $B=\sum \|a_m\|_2<\infty$ to conclude that 
 $\U a \in H^1_d$ with norm not exceeding $B$.  By \eqref{eq:bcmr},  we have  
 $$
\sup_{j\ge 0}\sup_{k\, ; \,   Q_{j,k} \subset Q_{0,m}}\bigg\{\sum_{\ell\, ;\, Q_{j,\ell} \subset Q_{0,0}} |\alpha^{j}_{k,\ell}|\bigg\} \le C(1+ |m|)^{-(n+s)}$$
and similarly exchanging the roles of $k$ and $\ell$. 
 Using Cauchy-Schwarz inequality,
\begin{align*}
 \|a_m\|^2_2& \le    \sum_{j=0}^\infty \sum_{k\, ;\, Q_{j,k} \subset Q_{0,m}}\bigg\{ \sum_{\ell\, ;\, Q_{j,\ell} \subset Q_{0,0}} |\langle a, \psi_{j,\ell}\rangle|^2 |\alpha^{j}_{k,\ell}| \bigg\}\bigg\{  \sum_{\ell\, ;\, Q_{j,\ell} \subset Q_{0,0}} |\alpha^{j}_{k,\ell}| \bigg\}   \\ 
 &   \le  C (1+ |m|)^{-(n+s)}   \sum_{j=0}^\infty \sum_{\ell\, ;\, Q_{j,\ell} \subset Q_{0,0}}  |\langle a, \psi_{j,\ell}\rangle|^2  \sum_{k\, ;\, Q_{j,k} \subset Q_{0,m}}   |\alpha^{j}_{k,\ell}|    
 \\
 &
\le C^2 (1+ |m|)^{-2(n+s)}   \sum_{j=0}^\infty \sum_{\ell\, ;\, Q_{j,\ell} \subset Q_{0,0}}  |\langle a, \psi_{j,\ell}\rangle|^2 
 \end{align*}
 and we are done provided one has a definition of  $\U$ on $H^1_d$. Let $\U_{J,N}$ be a partial sum obtained truncating  the sum defining $\U$ with $|j| \le J$ and $|k-\ell| \le 2^N$. It is immediate to define the action of $\U_{J,N}$ on all  $ H^1_d$ and we have
 $\|\U_{J, N} f\|_{H^1_d} \le C \| f\|_{H^1_d}$ for all $f \in H^1_d$ thanks to the previous calculations with $C$ independent of $J,N$ and $f$. Next,  by tedious but not difficult calculations refining the above estimates, one shows that 
 $\|\U_{J, N} f - \U f\|_{H^1_d} \to 0$ whenever $f$ is a finite linear combination of  $L^2$ dyadic atoms
 as $J, N \to \infty$.  Thus, we obtain the boundedness of $\U$ on a dense subspace of $H^1_d$   and we conclude by a density argument.

We next concentrate on $\W$. Once this is done, $\V$ is handled by observing that $\V^*$ is of the same type as $\W$. 
Recall that 
\begin{equation}
 \langle g, \W f \rangle =   \sum_{j=-\infty}^{+\infty}\sum_{k,\ell\in \ZZ^n}  \langle  g,  \phi_{j,k} \rangle  \gamma^{j}_{k,\ell}   \langle   \psi_{j,\ell},  f \rangle
 \end{equation}
 with $\gamma^j_{k,\ell}= c^j_{k,\ell}$ if $k\ne \ell$ and  $\gamma^j_{\ell,\ell}=- \sum_{k\ne \ell} c^j_{k,\ell}$. 
We decompose further $\W$ as 
$$\W=\sum_{R\in \NN^*} \W_R$$ where 
 \begin{equation}
 \langle g, \W_R f \rangle =   \sum_{j=-\infty}^{+\infty}\sum_{k,\ell\in \ZZ^n}  \langle  g,  \phi_{j,k} \rangle  \gamma^{j,R}_{k,\ell}   \langle   \psi_{j,\ell},  f \rangle
 \end{equation}
 \begin{equation}
\gamma^{j,R}_{k,\ell} =  \begin{cases}  c^j_{k,\ell}\ , & \mathrm{if}\   2^{R-1}\le |k-\ell| < 2^R ,
\\ - \sum_{2^{R-1} \le |m|<2^R} c^j_{k+m,k}\ ,  & \mathrm{if}\  k=\ell, 
\\ 0\ , & \mathrm{otherwise}.
\end{cases}
\end{equation}
Here,  for $x,y\in \RR^n$, $|x-y|= \sup(|x_1-y_1], \ldots, |x_n-y_n|)$. Let 
$$
\Gamma(R)= \sup_{j, k} \bigg\{\sum_{\ell} |\gamma^{j,R}_{k,\ell}| + |\gamma^{j,R}_{\ell,k}|\bigg\}.
$$
We notice that under \eqref{eq:bcmr}, we have $\Gamma(R)=O(2^{-Rs})$. 

\begin{lemma} 
For $R\geq 1$, we have: 
\begin{equation}
\|\W_{R}\|_{L^{2}\rightarrow L^{2}}\leq
CR^{\frac{1}{2}}\Gamma(R). 
\end{equation}
\begin{equation}
\|\W_{R}\|_{H^{1}_{d}\rightarrow L^{1}}  \leq
CR\,\Gamma(R).
\end{equation}
\begin{equation}
\|\W_{R}\|_{L^\infty\rightarrow BMO_d}\leq
CR\,\Gamma(R). 
\end{equation}
Hence, for $1<p<\infty$,
$$\|\W_{R}\|_{L^{p}\rightarrow L^{p}}\leq
CR\,\Gamma(R).$$
 \end{lemma}

It is clear that Theorem \ref{theor:bcmr}  for $\W$ follows at once from this lemma.

Let us begin the proof of this lemma by proving the $L^2$ boundedness. 
Write $f =  \sum_{j}f_{j}$ with $f_j=Q_j f$. Then,  
$$\|\W_{R}f\|_{2}^{2} \leq
\sum\limits_{j,j'}|\langle \W_{R}f_{j},
\W_{R}f_{j'}\rangle|.$$
  
First,  for each $j$, by expanding $f_j$ on the $\psi_{j,\ell}$, $\ell\in \ZZ^n$, and $\W_Rf_j$ on the 
$\phi_{j,k}$, $k\in \ZZ^n$,     Schur's lemma yields
$$
\|\W_{R}f_{j}\|_{2} \le \Gamma(R) \|f_j\|_2.
$$
Thus, using Cauchy-Schwarz inequality and   $\|f\|^2=\sum_j \|f_j\|_2^2$,  we have 
\begin{equation}
\sum\limits_{|j-j'|\leq R+2}|\langle \W_{R}f_{j},
\W_{R}f_{j'}\rangle| \le (2R+5) \Gamma(R)^2\|f\|_2^2.
\end{equation}

 It remains to handle the sum where $|j-j'|>
R+2$. It is enough to assume $j-j'>R+2$ and to show that 
\begin{equation}\label{1}
|\langle \W_{R}f_{j},
\W_{R}f_{j'}\rangle| \le C\,  \Gamma(R)^2\,  2^{\frac{j'-j+R}{2}} \|f_j\|_2 \|f_j'\|_2.\end{equation}
By dyadic scale invariance, assume also $j=0$, hence $-j'>R+2$. 
We have 
\begin{align*}
\langle \W_{R}f_{0},
\W_{R}f_{j'}\rangle &=\Big \langle \sum_{k,\ell} \langle f, \psi_{0,\ell} \rangle {}{\gamma^{0,R}_{k,\ell}} \phi_{0, k},  \sum_{k',\ell'} \langle f, \psi_{j',\ell'} \rangle {}{\gamma^{j',R}_{k',\ell'}} \phi_{j', k'} \Big\rangle
\\
&= \sum_{k'}  \bigg\{   \sum_{\ell} \langle f, \psi_{0,\ell}   \rangle \Big \langle \sum_k {}{\gamma^{0,R}_{k,\ell}}   \phi_{0,k}, \phi_{j',k'}\Big \rangle \bigg\} \bigg\{   \sum_{\ell'}{}{ \langle f, \psi_{j',\ell'} \rangle} \gamma^{j',R}_{k',\ell'}\bigg\}.
\end{align*}
Now the support of $\phi_{j',k'}$ is the cube $Q_{j',k'}= 2^{-j'}(k'+[0,1)^n)$ and for fixed $\ell$, the support of $\gamma^{0,R}_{k,\ell}$ is in the set of $k\in \ZZ^n$ such that $|k-\ell| < 2^R$. Thus, if $d(\ell, Q_{j',k'})>2^{R+2}$ then $\sum_k \gamma^{0,R}_{k,\ell}   \phi_{0,k}=0$ identically on $Q_{j',k'}$. Next, if $\ell \in Q_{j',k'}$ and $d(\ell, \RR^n\setminus Q_{j',k'})> 2^{R+2}$ then all the $\phi_{0,k}$ concerned in the sum have support inside $Q_{j',k'}$. Thus
$$
 \Big \langle \sum_k {}{\gamma^{0,R}_{k,\ell}}   \phi_{0,k}, \phi_{j',k'}\Big \rangle 
 = 2^{\frac{nj'}{2}} \sum_{k\in \ZZ^n}{}{ \gamma^{0,R}_{k,\ell} }=0
 $$
 by construction of the $\gamma$'s. Thus, for the sum in $\ell$, we have contribution only for $\ell\in E_{j',k'}$ defined as the set of those $\ell \in \ZZ^n$ with $d(\ell,\partial Q_{j',k'} ) \le 2^{R+2}$  (here, $\partial Q$ is the boundary of the cube $Q$) and the sum in $k$ inside the brackets reduces to those $k\in Q_{j',k'}$. 
Hence, we have
\begin{equation*}
|\langle \W_{R}f_{0},
\W_{R}f_{j'}\rangle | \le \sum_{k'}  \bigg\{   \sum_{\ell \in E_{j',k'} } | \langle f, \psi_{0,\ell}   \rangle|  \Big (\sum_{k\in Q_{j',k'}} |\gamma^{0,R}_{k,\ell}|   2^{\frac{nj'}{2}} \Big) \bigg\} \bigg\{   \sum_{\ell'}{| \langle f, \psi_{j',\ell'} \rangle | | \gamma^{j',R}_{k',\ell'}|}\bigg\}.
\end{equation*}
By Cauchy-Schwarz inequality, it suffices to estimate
$$
I= \bigg(\sum_{k'} \Big\{   \sum_{\ell'}{| \langle f, \psi_{j',\ell'} \rangle | | \gamma^{j',R}_{k',\ell'}|}\Big\}^2\ \bigg)^{1/2}
$$
and
$$
II= \bigg(\sum_{k'}  \Big\{   \sum_{\ell \in E_{j',k'} } | \langle f, \psi_{0,\ell}   \rangle|  \Big (\sum_{k\in Q_{j',k'}} |\gamma^{0,R}_{k,\ell}|   2^{\frac{nj'}{2}} \Big) \Big\}^2\ \bigg)^{1/2}.
$$
By Schur's lemma, we have
$$
I \le \Gamma(R) \bigg( \sum_{\ell'}  | \langle f, \psi_{j',\ell'} \rangle |^2  \bigg)^{1/2} = \Gamma(R) \|f_{j'}\|_2.
$$
Next,
\begin{align*}
II&=2^{\frac{nj'}{2}}   \bigg(\sum_{k'}  \bigg\{  \sum_{k\in Q_{j',k'}}  \sum_{\ell \in E_{j',k'} } | \langle f, \psi_{0,\ell}   \rangle|   |\gamma^{0,R}_{k,\ell}|  \bigg\}^2\ \bigg)^{1/2}\\
&\le  2^{\frac{nj'}{2}}   \bigg(\sum_{k'} \Big\{ \sum_{k\in Q_{j',k'}} \sum_{\ell \in E_{j',k'} } |\gamma^{0,R}_{k,\ell}|\Big\} \Big\{\sum_{k\in Q_{j',k'}}\sum_{\ell \in E_{j',k'} } | \langle f, \psi_{0,\ell}   \rangle|^2   |\gamma^{0,R}_{k,\ell}|   \Big\}\bigg)^{1/2}
\end{align*}
But, for fixed $k'$, since the cardinal of $E_{j',k'}$ is $O(2^{-j'(n-1)+R})$,
$$
\sum_{k\in Q_{j',k'}} \sum_{\ell \in E_{j',k'} } |\gamma^{0,R}_{k,\ell}| \le C\, \Gamma(R) \, 2^{-j'(n-1)+R}.
$$
Also
\begin{align*}
\sum_{k'} \sum_{k\in Q_{j',k'}}\sum_{\ell \in E_{j',k'} } | \langle f, \psi_{0,\ell}   \rangle|^2   |\gamma^{0,R}_{k,\ell}|  &\le \Gamma(R) \sum_{\ell \in \ZZ^n} | \langle f, \psi_{0,\ell}   \rangle|^2 
\Big\{ \sum_{k' \in \ZZ^n} {\bf 1}_{E_{j',k'}} (\ell)\Big\}
\\
&
\le  \Gamma(R) \sum_{\ell \in \ZZ^n} | \langle f, \psi_{0,\ell}   \rangle|^2 \ 2^n
\\
&
= 2^n\, \Gamma(R) \|f_0\|_2^2.
\end{align*}
All together
$$
II \le C\,  \Gamma(R)\,  2^{\frac{j'+R}2}\, \|f_0\|_2.
 $$
 and \eqref{1} is proved.

 \
 
 Next, we prove that $\W_R$ is bounded from $H^1_d$ to $L^1$ with norm $O(R\, \Gamma(R))$. To do that, we pick an $L^\infty$ dyadic atom $a$: it is supported in a dyadic cube $Q$, is bounded by $1/|Q|$ and is of mean 0. By scale and translation invariance, it suffices to assume that $Q=Q_{0,0}$.  Write $$a=\sum_{j'=0}^\infty\sum_{\ell\, ; \, Q_{j',\ell'} \subset Q_{0,0}} \langle a, \psi_{j',\ell'}\rangle  \psi_{j',\ell'}$$ and set
  $$a_1= \sum_{j'=R+1}^\infty\sum_{\ell'\, ; \, Q_{j',\ell'} \subset Q_{0,0}} \langle a, \psi_{j',\ell'}\rangle  \psi_{j',\ell'}.$$ Observe that $\|a_1\|_2\le \|a\|_2 \le \|a\|_\infty\le 1$. We have
 $$
 \W_R\,a_1= \sum_{j=-\infty}^\infty \sum_{k,\ell}   \langle a_1, \psi_{j,\ell} \rangle {}{\gamma^{j,R}_{k,\ell} } \phi_{j, k}  =\sum_{j=R+1}^\infty \sum_{k}\sum_{\ell; Q_{j,\ell} \subset Q_{0,0}} \langle a, \psi_{j,\ell} \rangle {}{\gamma^{j,R}_{k,\ell}} \phi_{j, k}
 $$
 and because $Q_{j,\ell} \subset [0,1]^n$, $|k-\ell|<2^R$ and $j\ge R+1$,  we have $Q_{j,k} \subset [-1,2]^n$ for all $(j,k)$ in the summation. Hence $\W_R\, a_1$ is supported in $[-1,2]^n$. Thus, the 
 $L^2$ estimate yields
 $$
 \|\W_R\,a_1\|_1\le C \|\W_R\, a_1\|_2\le CR^{1/2}\, \Gamma(R) \|a_1\|_2 \le CR^{1/2}\, \Gamma(R).
 $$
 Set $a_2=a-a_1$. Then, a straightforward estimate yields
 $$
 \|\W_R\,a_2\|_1 \le \sum_{j=0}^R  \sum_{k} \sum_{\ell\, ;\, Q_{j,\ell} \subset Q_{0,0}} |\gamma^{j,R}_{k,\ell}| 2^{-nj}\|a\|_\infty \le (R+1)\, \Gamma(R)
 $$
 by summing first in $k$, then in $\ell$ and in $j$. A truncation procedure with respect to the sum over  $j$ as for $\U$ allows to fully justify the boundedness of $\W_R$ from $H^1_d$ to $L^1$. We skip details which are easy. 
 
 \

Our last task is  prove that  $\W_R$ is bounded from $L^\infty$ to BMO$_d$  with norm $O(R\, \Gamma(R))$. Modulo a truncation procedure as above which is left to the reader, it suffices to show that $\W_R^*$ is bounded from $H^1_d$ to $L^1$.  So we pick again an $L^\infty$ dyadic atom $a$ and assume that it is supported in $Q=Q_{0,0}$. We
have 
 $$
\W_R^*\,a=  \sum_{j=-\infty}^{\infty} \sum_{k,\ell\in \ZZ^n} \langle a, \phi_{j,k}\rangle {\gamma^{j,R}_{k,\ell}} \psi_{j, \ell}= \sum_{j=0}^\infty \sum_{k\, ; \, Q_{j,k} \subset Q_{0,0}} \sum_{\ell}\langle a, \phi_{j,k}\rangle {\gamma^{j,R}_{k,\ell}} \psi_{j, \ell}$$
where we used that $a$ has support in $Q_{0,0}$ and mean 0.  We split the sum as $b_1+b_2$ according to $j\ge R+1$ or $j\le R$.  In the first case, we have as before, that 
 $Q_{j,k} \subset Q_{0,0}$, $j\ge R+1$ and $|k-\ell| <2^R$ imply that $Q_{j,\ell} \subset [-1,2]^n$ for all $(j,\ell)$ concerned by the summation. Also  $b_1$ can be written as $\widetilde \V_R (a)$ where $\widetilde \V_R$ is an operator of the same type as $\V_R$ with ``truncated'' coefficients (note that the $L^2$ bounds depends on a size estimate of the coefficients and on the nullity of the sum of the coefficients with respect to $\ell$ with $j$ and $k$ fixed). Thus, it is bounded on $L^2$ with bound $O(R^{1/2}\, \Gamma(R))$. Hence
$$
 \|b_1\|_1\le C \|b_1\|_2\le CR^{1/2}\, \Gamma(R) \|a\|_2 \le CR^{1/2}\, \Gamma(R).
 $$
For the $b_2$ part,  a straightforward estimate yields  a  bound
\begin{align*}
\| b_2 \|_1 & \le \sum_{j=0}^{R} \sum_{k\, ; \, Q_{j,k} \subset Q_{0,0}} \sum_{\ell\in \ZZ^n} \|a\|_\infty  2^{-nj} |\gamma^{j,R}_{\ell,k}|\\
& \le   (R+1) \, \Gamma(R)
\end{align*}
where $\Gamma (R)$ occurs by taking the sum in $\ell $ first. 

\begin{remark} It can be shown that $\W_R$ is bounded on BMO$_d$  with bound $O(R 2^{-Rs})$. Also $\W_R$ is  bounded from $H^1_d$ into $H^1$, the Hardy space on $\RR^n$ with a similar bound.   The proofs are a little more involved. However, it may not be bounded on $H^1_d$. The counterexample is the following: if $n=1$, set
$$
\langle g, \W f \rangle = (\langle g, \phi_{0,0}  \rangle - \langle g, \phi_{0,-1}  \rangle) \langle \psi_{0,0} , f \rangle 
$$
Then, observe that $\W_1=\W$ and $\W_1(\psi_{0,0})(x)= \phi_{0,0}(x)  - \phi_{0,-1}(x)= \phi(x) - \phi(x+1) \notin H^1_d$ since it does not vanish on $\RR^+$ and $\RR^-$ which is necessary. 
\end{remark}

\end{document}